\documentclass[3p,times]{elsarticle}
\usepackage{amssymb, color}
\usepackage{amsmath}
 \usepackage{pifont}
 \usepackage[utf8]{inputenc}
\usepackage[frenchb]{babel}
\usepackage[colorlinks,,breaklinks]{hyperref}
\journal{arXiv}
\makeatletter
\def\LaTeX{\leavevmode L\raise.42ex
    \hbox{\kern-.3em\size{\sf@size}{0pt}\selectfont A}\kern-.15em\TeX}
\makeatother

\newcommand{\BibTeX}{{\rm B\kern-.05em{\sc
          i\kern-.025emb}\kern-.08em\TeX}}

\makeatletter
\def\@currentlabel{2.1}\label{e:dispaa}
\def\@currentlabel{2.21}\label{e:dispau}
\def\@currentlabel{2.22}\label{e:dispav}
\def\@currentlabel{2.23}\label{e:dispaw}
\def\@currentlabel{2.24}\label{e:dispax}
\def\theequation{\thesection.\@arabic\c@equation}
\makeatother

\renewcommand{\theequation}{\arabic{section}.\arabic{equation}}

\newcommand{\R}{\mathbb R}

\def \D{\Delta}

\newtheorem{thm}{Theorem} [section]

\newtheorem{lem}{Lemma} [section]
\newtheorem{prop}{Proposition} [section]

\newtheorem{definition}{Definition} [section]

\newtheorem{rem}{Remark}[section]

\renewcommand{\theequation}{\thesection.\arabic{equation}}
\renewcommand{\thesection}{\arabic{section}}
\renewcommand{\theequation}{\thesection.\arabic{equation}}
\let\ssection=\section\renewcommand{\section}{\setcounter{equation}{0}\ssection}

\begin{document}
\begin{large}
\begin{frontmatter}

\title{Liouville-type theorems with finite Morse index for $\D_{\lambda}$-Laplace operator}

\author[br]{Belgacem Rahal }
\ead{rahhalbelgacem@gmail.com}

\address[br]{Facult\'e des Sciences, D\'epartement de Math\'ematiques, B.P 1171 Sfax 3000, Universit\'e de Sfax, Tunisia.}

\begin{abstract} In this paper we study solutions, possibly unbounded and sign-changing, of the following problem 
  $$-\D_{\lambda} u=|x|_{\lambda}^a |u|^{p-1}u,\quad \mbox{in}\,\, \R^n,\;n\geq 1,\; p>1,\; \mbox{ and }\; a \geq 0,$$
where $\D_{\lambda}$ is a strongly degenerate elliptic operator, the functions $\lambda=(\lambda_1, ..., \lambda_k) : \; \R^n \rightarrow \R^k$ satisfies some certain conditions, and $|.|_{\lambda}$  the homogeneous norm associated to the $\D_{\lambda}$-Laplacian.
 We prove various Liouville-type theorems for smooth solutions under the assumption that they are
stable or stable outside a compact set of $\R^n$. First, we establish the standard integral
estimates via stability property to derive the nonexistence results  for stable solutions. Next, by mean of the Pohozaev
identity, we deduce the Liouville-type theorem for solutions stable outside a compact set.
\end{abstract}
\begin{keyword}
\begin{large}
Liouville-type theorems\sep $\D_{\lambda}$-Laplace operator \sep Stable solutions \sep Stability outside a compact set \sep Pohozaev
identity.
\PACS {Primary: 35J55, 35J65; Secondary: 35B65.}
\end{large}
\end{keyword}
\end{frontmatter}
\section{\large{Introduction and main results}}
The Liouville type theorem is the nonexistence of solutions in the entire space or in half-space. The classical Liouville type theorem stated that a bounded harmonic (or holomorphic) function defined in entire space must be constant. This theorem, known as Liouville theorem, was first announced in 1844 by Liouville \cite{Liouville} for the special case of a doubly-periodic function. Later in the same year, Cauchy \cite{Cauchy} published the first proof of the above stated theorem. This classical result has been extended to nonnegative solutions of the semilinear elliptic equation
\begin{align}\label{equation0}
  -\Delta u=|u|^{p-1}u\ \ \   \ \ \mbox{in}\,\, \R^n,\;  p>1,
\end{align}
in the whole space $\R^n$ by Gidas and Spruck \cite{Gidas1, Gidas2} see also the paper of Chen and Li \cite{Chen}. They proved that if $1<p< \frac{n+2}{n-2}$, then the above equation only has the trivial solution $u \equiv 0$ and this result is optimal. In an elegant paper, Farina \cite{Farina} proved that nontrivial finite Morse index solutions (whether positive or
sign changing) to \eqref{equation0} exists if and only if $p\geq p_{c}(n)$ and $n\geq 11$, or $ p=\frac{n+2}{n-2}$ and $n\geq 3$, where $ p_{c}(n) $ is the so-called Joseph-Lundgren exponent. The study of stable solutions in the H\'{e}non type elliptic equation: $-\Delta u= |x|^a |u|^{p-1}u,\;\; \mbox{in}\,\, \R^n,\;  p>1\; \mbox{and}\; a> -2$ has been studied recently,
Wang and Ye \cite{WY} gave a complete classification of stable weak solutions and those of finite Morse index solutions.

In the past years, the Liouville property has been refined considerably and emerged as one of the most powerful tools in the study of initial and boundary value problems for nonlinear PDEs. It turns out that one can obtain from Liouville-type theorems a variety of results on qualitative properties of solutions such as universal, pointwise, a priori estimates of local solutions; universal and singularity estimates; decay estimates; blow-up rate of solutions of nonstationary problems, etc., see \cite{Polacik, Quittner} and references therein.

Liouville-type theorems for degenerate elliptic equations have been attracted the interest of many mathematicians. The classical Liouville theorem was generalized to $p$-harmonic functions on the whole space $\R^n$ and on exterior domains by Serrin and Zou \cite{Serrin}, see also \cite{Cuccu} for related results. The Liouville theorems for some linear degenerate elliptic operators such as $X$-elliptic operators, Kohn-Laplacian (and more general sublaplacian on Carnot groups) and degenerate Ornstein-Uhlenbeck operators were proved in \cite{Kogoj1, Kogoj2}.


 More recently, Yu \cite{Yu} studied the equation
$$- L_{\alpha} u = f(u)\;\;\mbox{ in }\;\; \R^{n_1} \times \R^{n_2}, $$
where $L_{\alpha}= \D_{x} + (1+\alpha)^2 \D_{y}$, $\alpha >0$ and $Q= n_1 + (1+\alpha ) n_2$ is the \textit{homogeneous dimension} of the space. Under some assumptions on the nonlinear term $f$ , he showed that the above equation possesses no positive solutions and the main technique used is the moving plane method in the integral form.

In this paper, we are concerned with the Liouville-type theorems for the following problem \begin{align}\label{equation}
  -\Delta_{\lambda} u=|x|^{a}_{\lambda} |u|^{p-1}u,\quad \mbox{in}\,\, \R^n:= \R^{n_1}\times \R^{n_2}\times ...\times \R^{n_k},
\end{align}
where $n\geq 1$, $a \geq 0$, $p>1$, $$ \D_{\lambda} = \lambda^2_{1} \D_{x^{(1)}} +...+ \lambda^2_{k} \D_{x^{(k)}},\quad |x|_{\lambda} := \left(\sum_{j=1}^k \prod_{i\neq j} \lambda_i^2(x) \epsilon^2_j |x^{(j)}|^2 \right)^{\frac{1}{2\sigma }},$$ $\sigma= 1+\sum_{i=1}^k (\epsilon_i-1)$, $1\leq \epsilon_1 \leq ...\leq \epsilon_k$, $x =(x^{(1)}, ..., x^{(k)}) \in \R^n$.
Here the functions $\lambda_i: \R^n \rightarrow \R$ are continuous, strictly
positive and of class $C^1$ outside the coordinate hyperplanes, i.e. $\lambda_i >0$, $i=1,...,k$ in $\R^n \backslash \prod$, where $\prod = \{ x=(x_1,...,x_n)\in \R^n : \prod_{i=1}^n x_i=0\}$, and $\D_{x^{(i)}}$ denotes the classical
Laplacian in $\R^{n_i}$, $i=1,...,k$. As in \cite{Kogoj3} we assume that $\lambda_i$ satisfy the following properties:\\\\
\textbf{$(H_1)$} $\lambda_1(x)=1$, $\lambda_i(x)= \lambda_i(x^{(1)}, ..., x^{(i-1)})$, $i=2,..., k$.\\\\
\textbf{$(H_2)$} For every $x\in \R^n$, $\lambda_i(x)= \lambda_i(x^*)$, $i=1, ..., k$, where $x^*=( |x^{(1)}|,...,|x^{(k)}|)$ if $x=(x^{(1)},..., x^{(k)})$. \\\\
\textbf{$(H_3)$} There exists a group of dilations $\{\delta_t\}_{t>0}$ $$\delta_t : \R^n \rightarrow \R^n , \; \delta_t(x)= \delta_t(x^{(1)},..., x^{(k)})= (t^{\epsilon_1}x^{(1)},..., t^{\epsilon_k}x^{(k)}), $$
where $1 \leq \epsilon_1 \leq \epsilon_2 \leq ...\leq \epsilon_k$, such that $\lambda_i$ is $\delta_t$-homogeneous of degree $\epsilon_i-1$, i.e. $$\lambda_i (\delta_t(x))=  t^{\epsilon_i -1} \lambda_i(x),\;\; \forall \; x\in \R^n , \; t>0, \; i=1,..., k. $$
This implies that the operator $\D_{\lambda}$ is $\delta_t$-homogeneous of degree two, i.e. $$\D_{\lambda}(u(\delta_t(x)))= t^2 (\D_{\lambda} u) (\delta_t(x)),\; \; \forall \; u \in C^{\infty}(\R^n) .$$
We denote by $Q$ the \textit{homogeneous dimension } of $\R^n$ with respect to the group of dilations $\{\delta_t\}_{t>0}$, i.e. $$Q:= \epsilon_1 n_1+ \epsilon_2 n_2 +...+\epsilon_k n_k.$$
The $\D_{\lambda}$-Laplace operator was first introduced by Franchi and Lanconelli \cite{Franchi}, and recently reconsidered in \cite{Kogoj3} under an additional assumption that the operator is homogeneous of degree two with respect to a group dilation in $\R^n$. It was proved in \cite{AnhMy}, that the autonomous case, i.e. $a=0$, \eqref{equation} has no positive classical solution if $1< p\leq \frac{Q}{Q-2}$, with $Q= \epsilon_1 + \epsilon_2  +...+\epsilon_n $, ($n_i=1$, $i=1, ..., n$).

 The $\D_{\lambda}$-operator contains many degenerate elliptic operators. We now give some examples of $\D_{\lambda}$-Laplace operators (see also \cite{Kogoj3}). We use the following notation: we split $\R^n$ as follows $\R^n=\R^{n_1} \times ...\times \R^{n_k}$ and write $$x=(x^{(1)}, ..., x^{(k)}), \; x^{(i)}=(x^{(i)}_1, ..., x^{(i)}_{n_i}) \in \R^{n_i},$$
$$|x^{(i)}|^{2}= \sum_{j=1}^{n_i} |x^{(i)}_j|^2,\;\;\; i=1,2,...,k.$$
We denote the classical Laplace operator in $\in \R^{n_i}$ by $$\D_{x^{(i)}} = \sum_{j=1}^{n_i} \partial^2_{x^{(i)}_j}.$$
\textbf{Example 1.} Let $\alpha$ be a real positive constant and $k = 2$. We consider the Grushin-type
operator $$ \D_{\lambda}= \D_x + |x|^{2 \alpha } \D_y,$$
where $\lambda = (\lambda_1,\lambda_2)$ with $$\lambda_1(x)=1, \quad \lambda_2(x)= |x^{(1)}|^{\alpha },\quad x \in \R^{n_1} \times \R^{n_2}.$$
Our group of dilations is
$$\delta_t(x)= \delta_t(x^{(1)},x^{(2)})= (tx^{(1)}, t^{\alpha+1}x^{(2)}),$$
and the homogenous dimension with respect to $(\delta_t)_{t>0}$ is $Q = n_1 + (\alpha + 1)n_2$.\\\\
\textbf{Example 2.}
Given a multi-index $\alpha=(\alpha_1,..., \alpha_{k-1}) $, $\alpha_j \geq 1 $, $j=1,...,k-1$, define $$\D_{\alpha}:= \D_{x^{(1)}}+ |x^{(1)}|^{2 \alpha_1} \D_{x^{(2)}} + ...+ |x^{(k-1)}|^{2 \alpha_{k-1}} \D_{x^{(k)}}.$$
Then $\D_{\alpha}= \D_{\lambda}$ with $\lambda= (\lambda_1, ..., \lambda_k )$ and $\lambda_i= |x^{(i-1)}|^{\alpha_{i-1}} $, $i=1,...,k$. Here we agree to let $|x^{(0)}|^{\alpha_{0}}=1$. A group of dilations for which $\lambda$ satisfies $(H_3)$ is given by
$$\delta_t : \R^n \rightarrow \R^n , \; \delta_t(x)= \delta_t(x^{(1)},..., x^{(k)})= (t^{\epsilon_1}x^{(1)},..., t^{\epsilon_k}x^{(k)}), $$ with $\epsilon_1=1$ and $\epsilon_i= \alpha_{i-1} \epsilon_{i-1}+1$, $i=2,...,k$. In particular, if $\alpha_1=...=\alpha_{k-1}=1$, the operator $\D_{\alpha}$ and the dilation $\delta_t$ becomes, respectively
$$\D_{\alpha}= \D_{x^{(1)}}+ |x^{(1)}|^{2 } \D_{x^{(2)}} + ...+ |x^{(k-1)}|^{2 } \D_{x^{(k)}},$$
and
$$\delta_t(x)= (t x^{(1)},t^{2}x^{(2)},..., t^{k}x^{(k)}).$$
\textbf{Example 3.} Let $\alpha$, $\beta$  and  $\gamma$ be positive real constants. For the operator
 $$\D_{\lambda} = \D_{x^{(1)}} + |x^{(1)}|^{2\alpha } \D_{x^{(2)}}+ |x^{(1)}|^{2\beta } |x^{(2)}|^{2\gamma } \D_{x^{(3)}},$$
 where $\lambda = (\lambda_1,\lambda_2,\lambda_3)$ with $$\lambda_1(x)=1, \quad \lambda_2(x)= |x^{(1)}|^{\alpha }     , \quad \lambda_3(x)=|x^{(1)}|^{\beta } |x^{(2)}|^{\gamma },\quad x \in \R^{n_1} \times \R^{n_2} \times\R^{n_3},$$
 we find the group of dilations
 $$\delta_t(x)= \delta_t(x^{(1)},x^{(2)}, x^{(3)})= (tx^{(1)}, t^{\alpha+1}x^{(2)} , t^{\beta + (\alpha+1)\gamma +1}x^{(3)}).$$

The aim of the present paper was to establish the Liouville-type theorems with finite Morse index for the equation \eqref{equation}. In order to state our results we need the following:
\begin{definition}
 We say that a solution $u$ of \eqref{equation} belonging to $C^2(\R^n)$\\
 $\bullet$ is stable, if
$$Q_u(\psi):= \int_{\R^n} |\nabla_{\lambda} \psi|^2 - p\int_{\R^n}  |x|^{a}_{\lambda} |u|^{p-1}\psi^2 \geq0\;,\;\;\forall\; \psi\in  C_c^1( \R^n),$$
where $\nabla_{\lambda}= (\lambda_1 \nabla_{x^{(1)}},..., \lambda_k \nabla_{x^{(k)}})$.\\
$\bullet$ has Morse index equal to $K \geq 1$ if $K$ is the maximal dimension of a subspace $X_K$ of $C^1_c(\R^n)$ such that $ Q_u (\psi)< 0$ for any $\psi \in X_K\backslash \{0\}$.\\
$\bullet$ is stable outside a compact set $\mathcal{K}\subset \R^n$ if $ Q_u (\psi)\geq0$ for any $\psi \in C^1_c(\R^n\backslash \mathcal{K})$.
\end{definition}
\begin{rem}
\textbf{a)} Clearly, a solution stable if and only if its Morse index is equal to zero.\\
\textbf{b)} It is well know that any finite Morse index solution $u$ is stable outside a compact set $\mathcal{K}\subset \R^n$. Indeed, there exists $m_0\geq1$ and $X_{m_0}:= \mbox{Span}\{\phi_1,..., \phi_{m_0}\} \subset C^1_c(\R^n)$ such that $Q_u (\phi)<0$ for any $\phi \in X_{m_0} \backslash\{0\}$. Hence, $Q_u (\psi)\geq0$ for every $\psi \in C^1_c(\R^n \backslash \mathcal{K}) $, where $\mathcal{K}:= \displaystyle{ \cup_{j=1}^{m_0}} supp(\phi_{j})$.
\end{rem}

In the following, we state Liouville-type results for solutions $u \in C^2(\R^n)$  of \eqref{equation}. In what follows,
we divide our study to stable solutions and solutions which are stable outside a compact set.

\subsection{\textbf{\large{Stable solutions}}}

 To state the following result we need to introduce some notation. We set $\Gamma_M(p)=2p-1+2\sqrt{p(p-1)}$ and denote by $\Omega_{R}= B_1(0,  R^{\epsilon_1})\times B_2(0,  R^{\epsilon_2}) \times ...\times B_k(0, R^{\epsilon_k})$, where $B_i(0,  R^{\epsilon_i})\subset \R^{n_i}$, $i=1,...,k$,  the balls of center $0$ and radius $R^{\epsilon_i}$.
\begin{prop}\label{proposition1}
Let $u \in C^2(\R^n)$ be a stable solution of \eqref{equation}. Then, for any $\gamma \in \left[1, \Gamma_M(p)\right)$, there exists a positive constant $C$ independent of $R$, such that
\begin{eqnarray} \label{stable} \int_{\Omega_R} \left(|x|^{a}_{\lambda} |u|^{p+\gamma}+ |\nabla_{\lambda} ( |u|^{\frac{\gamma-1}{2}}u)|^2 \right) dx \leq C R^{Q-  \frac{2(p+\gamma)+ (\gamma+1)a }{p-1}},\quad \mbox{for all} \; R>0.
\end{eqnarray}
\end{prop}
Proposition \ref{proposition1} provides an important estimate on the integrability of $u$ and $\nabla_\lambda u$. As we will see, our nonexistence results will follow by showing that the right-hand side of \eqref{stable} vanishes under the right assumptions on $p$ when $R\rightarrow +\infty$. More precisely, as a corollary of Proposition \ref{proposition1}, we can state our first Liouville type theorem.
\begin{thm}\label{th1}
Let $u \in C^2(\R^n)$ be a stable solution of \eqref{equation} with, \begin{eqnarray*}
 p_c(Q,a)=
\begin{cases}
+\infty\;\;&\text{if  $Q\leq 10 + 4 a $},\\
 \frac{(Q-2)^2-2(a+2)(a+Q)+2 \sqrt{(a+2)^3(a+2Q-2)}}{(Q-2)(Q-4a-10)}\;\; &\text{if  $Q > 10 + 4 a $}.
\end{cases}
\end{eqnarray*} Then $u\equiv 0$.
\end{thm}
\subsection{\textbf{\large{Solutions which are stable outside a compact set}}}

In this subsection we prove some integral identities extending to the $\D_{\lambda}$ setting the classical Pohozaev identity for semilinear
Poisson equation \cite{Pohozaev}. Pohozaev identity has been extended by several authors to general elliptic equations and systems,
both in Riemannian and sub-Riemannian context, see, e.g., \cite{Bozhkov, Garofalo, Pucci} and the references therein. To prove our identities we
closely follow the original procedure of Pohozaev, just replacing the vector field $ P= \sum_{i=1}^n x_i \partial_{x_i} $ in \cite{Pohozaev}, page $1410$], by $$ T= \sum_{i=1}^k \epsilon_i x^{(i)} \nabla_{x^{(i)}},$$
the generator of the group of dilation $(\delta_t)_{t\geq 0}$ in \textbf{$(H_3)$}(we say that $T$ generates $(\delta_t)_{t\geq 0}$ since a function $u$ is $\delta_t$-homogeneous of degree $m$ if and only if $T u = m u$).
\begin{prop}\label{proposition2}
Let $u \in C^2(\R^n)$ be a solution of \eqref{equation} and $\phi \in C^1_c(\Omega_R)$. If $T\left( |x|_{\lambda}\right)= |x|_{\lambda}$, then
\begin{multline} \label{pohozaev}
 \int_{\Omega_R}  \left[ \frac{Q-2}{2}  |\nabla_{\lambda} u|^2 -\frac{Q+a}{p+1} |x|^{a}_{\lambda}|u|^{p+1} \right]\phi = \int_{\Omega_R} \left[\nabla_{\lambda} u \nabla_{\lambda} \phi  T(u) +\left[- \frac 12 |\nabla_{\lambda} u|^2 + \frac{|x|^{a}_{\lambda}}{p+1} |u|^{p+1}\right] T(\phi)\right].
\end{multline}
\end{prop}
Thanks to Proposition \ref{proposition2}, we derive
\begin{thm}\label{th2}
Let $u \in C^2(\R^n)$ be a solution of \eqref{equation} which is stable outside a compact set of $\R^n$, with \begin{eqnarray*}
p_s(Q,a)=\begin{cases}
+\infty  \; \; &\text{if \; $Q\leq 2$,}\\
\frac{Q+2+2a}{Q-2} \; \; &\text{if \; $Q> 2$.}
\end{cases}
\end{eqnarray*} If $T\left( |x|_{\lambda}\right)= |x|_{\lambda}$, then $u\equiv 0$.
\end{thm}

\section{ \large{Example which satisfies $T\left( |x|_{\lambda}\right)= |x|_{\lambda}$} }
The degenerate elliptic operators we consider are of the form $$ \D_{\lambda} = \lambda^2_{1} \D_{x^{(1)}} +...+ \lambda^2_{k} \D_{x^{(k)}}. $$
We denote by $|x^{(j)}|$ the euclidean norm of $x^{(j)} \in \R^{n_j}$ and assume the functions $\lambda_{i}$ are
of the form
\begin{eqnarray} \label{equal}
\lambda_{i}(x)= \prod_{j=1}^k |x^{(j)}|^{\alpha_{ij}},\quad i=1,...,k,
\end{eqnarray}
such that\\
1) $\alpha_{ij} \geq 0$ for $i = 2, . . . , k$, $j = 1, . . . , i-1$.\\
2) $\alpha_{ij} = 0$ for $j\geq i$.\\
3) $\sum_{l=1}^k \epsilon_l \alpha_{jl} = \epsilon_j-1$,  $j=1,...,k$ with $1=\epsilon_1\leq \epsilon_2 \leq ...\leq \epsilon_k$.\\
Clearly, $\lambda_i$ is $\delta_t$-homogeneous of degree $\epsilon_i-1$ with respect to a group of dilations $\{\delta_t\}_{t>0}$ $$\delta_t : \R^n \rightarrow \R^n , \; \delta_t(x)= \delta_t(x^{(1)},..., x^{(k)})= (t^{\epsilon_1}x^{(1)},..., t^{\epsilon_k}x^{(k)}). $$
Now, using the relation $\sum_{l=1}^k \epsilon_l \alpha_{jl} = \epsilon_j-1$, we get $T\left( |x|_{\lambda}\right)= |x|_{\lambda}$ is satisfied.\\

This paper is organized as follows. In section $3$, we give the proof of Proposition \ref{proposition1} and Theorem \ref{th1}. Section $4$ is devoted to the proof of Proposition \ref{proposition2} and Theorem \ref{th2}.
\section{\large{The Liouville theorem for stable solutions: proof of Theorem \ref{th1}}}
In this section we prove all the results concerning the classification of stable solutions, i.e., Proposition \ref{proposition1} and Theorem \ref{th1}. First, to prove Proposition \ref{proposition1}, we need the following technical Lemma.

Let $R>0$, $\Omega_{2R}= B_1(0, 2 R^{\epsilon_1})\times B_2(0, 2 R^{\epsilon_2}) \times ...\times B_k(0, 2 R^{\epsilon_k})$, where $B_i(0, 2 R^{\epsilon_i})\subset \R^{n_i}$, $i=1,...,k$, and consider $k$ functions $\psi_{1,R}$,..., $\psi_{k,R}$ such that $$\psi_{1,R}(r^{(1)})=\psi_{1}\left(\frac{r^{(1)}}{R^{\epsilon_1}}\right),\;...,\; \psi_{k,R}(r^{(k)})=\psi_{k}\left(\frac{r^{(k)}}{R^{\epsilon_k}}\right),$$
with $\psi_{1,R},\; ...\; \psi_{k,R} \in C^{\infty}_c([0, + \infty)),\; 0\leq \psi_{1,R},\;...\; \psi_{k,R} \leq 1,$ $$ \psi_i(t)=\begin{cases} 1\; & \text{in \;$[0,\;1]$},\\
0\; & \text{in \;$[2,\;+\infty )$},
\end{cases}$$
and for some constant $C>0$ and $\psi_{1,R},\; ...\; \psi_{k,R}$ satisfy
\begin{eqnarray*} \label{ineq1}
\left|\nabla_{x^{(1)}} \psi_{1,R}\right| \leq C R^{-\epsilon_1},\; ...,\; \left|\nabla_{x^{(k)}} \psi_{k,R}\right| \leq C R^{-\epsilon_k},
\end{eqnarray*}
\begin{eqnarray*} \label{ineq2}
\left|\D_{x^{(1)}} \psi_{1,R}\right| \leq C R^{-2\epsilon_1},\; ...,\; \left|\D_{x^{(k)}} \psi_{k,R}\right| \leq C R^{-2\epsilon_k},
\end{eqnarray*}
where $r^{(i)}= |x^{(i)}|$, $i=1,..., k$.

\begin{lem} \label{lem1}
\textbf{(1)} There exists a constant $C>0$ independent of $R$ such that

 \textbf{a)} $|\lambda_{i}(x)| \leq C R^{\epsilon_i-1},\; \forall \; x \in \Omega_{2R},\; i=1,...,k.$

\textbf{b)} $|\nabla_{\lambda} \psi_R|^2 + |\D_{\lambda} \psi_R|\leq C R^{-2}$, where
$\psi_R= \prod_{i=1}^k \psi_{i,R}$.\\\\
\textbf{(2)} The homogeneous norm, $|.|_{\lambda}$, is $\delta_t$-homogeneous of degree one, i.e. $$|\delta_t(x)|_{\lambda}=  t |x|_{\lambda},\;\; \forall \; x\in \R^n , \; t>0. $$
\textbf{(3)} There exists a constant $C>0$ independent of $R$ such that$$|x|_{\lambda}\leq C R,\;\forall \; x \in \Omega_{2R}.$$
\end{lem}
\textbf{Proof.}\\ \textit{Proof of (1)\;a).} For any $x=(x^{(1)},..., x^{(k)})\in \Omega_{2R}$, we have $x^{(i)} \in B_i(0, 2R^{\epsilon_i})$, $i=1,...,k$, this implies $\frac{|x^{(i)}|}{R^{\epsilon_i}}\leq 2$, $i=1,...,k$. Therefore, if we write $$x=(x^{(1)},..., x^{(k)})= \left( R^{\epsilon_1} \times \frac{x^{(1)}}{R^{\epsilon_1}}, ..., R^{\epsilon_k} \times \frac{x^{(k)}}{R^{\epsilon_k}}\right), $$
and let $y=(y^{(1)}, ..., y^{(k)})= \left( \frac{x^{(1)}}{R^{\epsilon_1}}, ..., \frac{x^{(k)}}{R^{\epsilon_k}}\right)$, then $y \in \overline{\Omega_{2}}$. Hence by assumption $(H_3)$ made on functions $\lambda_i$, we get
\begin{eqnarray} \label{p1}
\lambda_i(x)&=& \lambda_i(R^{\epsilon_1} y^{(1)}, ..., R^{\epsilon_k} y^{(k)} )\nonumber\\&=& R^{\epsilon_i-1} \lambda_i(y^{(1)}, ..., y^{(k)}) \nonumber\\&=& R^{\epsilon_i-1} \lambda_i(y).
\end{eqnarray}
Moreover, since $\lambda_i$, $i=1, ..., k$ are continuous, then
\begin{eqnarray} \label{p2}
|\lambda_i(y)|\leq C, \;\; \forall \; y \in \overline{\Omega_{2}}.
\end{eqnarray}
Therefore, from \eqref{p1} and \eqref{p2}, we obtain $$|\lambda_{i}(x)| \leq C R^{\epsilon_i-1},\; \forall \; x \in \Omega_{2R},\; i=1,...,k.$$
\textit{Proof of (1)\;b).} Using assumption $(H_2)$ made on functions $\lambda_i$, $i=1,...,k$, with $r=( r^{(1)}, ..., r^{(k)})=(|x^{(1)}|, ..., |x^{(k)}|)$, we have $$\lambda_1(r)=1,\; \lambda_i(r)=\lambda_i(r^{(1)}, ..., r^{(i-1)}),\; \forall\; i=2,...,k.$$
If we denote by $\psi_R= \prod_{i=1}^k \psi_{i,R}$, we get
\begin{eqnarray*}
\nabla_{\lambda} \psi_R &=&\left( \lambda_1(r) \nabla_{x^{(1)}} \psi_R,\; ...,\; \lambda_k(r) \nabla_{x^{(k)}}\psi_R\right)\nonumber\\&=&
\left(\lambda_1(r)\nabla_{x^{(1)}}\psi_{1,R}\prod_{i=2}^k \psi_{i,R}, ...,  \lambda_k(r)\nabla_{x^{(k)}}\psi_{k,R}\prod_{i=1}^{k-1} \psi_{i,R} \right),
\end{eqnarray*}
and
\begin{eqnarray*}
\D_{\lambda} \psi_R &=& \lambda^2_1(r) \D_{x^{(1)}} \psi_R +...+ \lambda^2_k(r) \D_{x^{(k)}} \psi_R  \nonumber\\ &=& \lambda^2_1(r) \D_{x^{(1)}}\psi_{1,R}\prod_{i=2}^k \psi_{i,R} +...+ \lambda^2_k(r) \D_{x^{(k)}}\psi_{k,R}\prod_{i=1}^{k-1} \psi_{i,R}.
\end{eqnarray*}
Since $|\lambda_{i}(r)|=|\lambda_{i}(x)| \leq C R^{\epsilon_i-1}$, $\forall \; x \in \Omega_{2R}$, $i=1,...,k$, then there exists a constant $C>0$ independent of $R$ such that
$$|\nabla_{\lambda} \psi_R|^2\leq C R^{-2}\; \; \mbox{and}\;\; |\D_{\lambda} \psi_R|\leq C R^{-2}.$$
\textit{Proof of (2).} Let $x \in \R^n$. The homogeneity of the functions $\lambda_i$ implies that
\begin{eqnarray}\label{nouveau}
|\delta_t(x)|_{\lambda} :&=& \left(\sum_{j=1}^k \prod_{i\neq j} (\lambda_i(\delta_t(x)))^2 \epsilon^2_j |t^{\epsilon_j} x^{(j)}|^2 \right)^{\frac{1}{2(1+\sum_{i=1}^k (\epsilon_i-1))}} \nonumber\\&=& \left(\sum_{j=1}^k \prod_{i\neq j} t^{2\epsilon_j} t^{2(\epsilon_i-1)} (\lambda_i(x))^2 \epsilon^2_j | x^{(j)}|^2 \right)^{\frac{1}{2(1+\sum_{i=1}^k (\epsilon_i-1))}} \nonumber\\&=& \left( t^{2(1+\sum_{i=1}^k(\epsilon_i-1))}\sum_{j=1}^k \prod_{i\neq j}  (\lambda_i(x))^2 \epsilon^2_j | x^{(j)}|^2 \right)^{\frac{1}{2(1+\sum_{i=1}^k (\epsilon_i-1))}}\nonumber\\&=& t |x|_{\lambda}
\end{eqnarray}
\textit{Proof of (3).} For any $x=(x^{(1)},..., x^{(k)})\in \Omega_{2R}$, we have $x^{(i)} \in B_i(0, 2R^{\epsilon_i})$, $i=1,...,k$, this implies $\frac{|x^{(i)}|}{R^{\epsilon_i}}\leq 2$, $i=1,...,k$. Therefore, if we write $$x=(x^{(1)},..., x^{(k)})= \left( R^{\epsilon_1} \times \frac{x^{(1)}}{R^{\epsilon_1}}, ..., R^{\epsilon_k} \times \frac{x^{(k)}}{R^{\epsilon_k}}\right), $$
and let $y=(y^{(1)}, ..., y^{(k)})= \left( \frac{x^{(1)}}{R^{\epsilon_1}}, ..., \frac{x^{(k)}}{R^{\epsilon_k}}\right)$, then $y \in \overline{\Omega_{2}(0)}$.\\
Using \eqref{nouveau}, we get
\begin{eqnarray*}
|x|_{\lambda}&=& |(R^{\epsilon_1} y^{(1)}, ..., R^{\epsilon_k} y^{(k)} )|_{\lambda}\nonumber\\&=& R |(y^{(1)}, ..., y^{(k)})|_{\lambda} \nonumber\\&=& R|y|_{\lambda}.
\end{eqnarray*}
Since $|\lambda_{i}(y)| \leq C$, $\forall \; y \in \overline{\Omega_{2}}$, $i=1,...,k$, then there exists a constant $C>0$ independent of $R$ such that $$|x|_{\lambda} \leq C R,\;\; \forall \; x \in  \Omega_{2R}.$$
This completes the proof of Lemma \ref{lem1}. \qed\\\\
\textbf{Proof of Proposition \ref{proposition1}.} The proof follows the main lines of the demonstration of proposition $4$ in \cite{Farina}, with more modifications. We split the proof into four steps: \\
\textit{\textbf{Step 1.}} For any $\phi \in C^2_c(\R^n)$ we have
\begin{eqnarray}\label{e1}
\int_{\R^n}   |\nabla_{\lambda} (|u|^{\frac{\gamma-1}{2}}u)|^2 \phi^2 dx&=& \frac{(\gamma +1)^2}{4 \gamma}\int_{\R^n} |x|^{a}_{\lambda} |u|^{p+\gamma} \phi^2dx + \frac{\gamma +1}{4 \gamma} \int_{\R^n} |u|^{\gamma+1} \Delta_{\lambda} (\phi^2)dx.
\end{eqnarray}
Multiply equation \eqref{equation} by $|u|^{\gamma-1} u \phi^2$
and integrate by parts to find
$$\gamma \int_{\R^n} |\nabla_{\lambda} u|^2 |u|^{\gamma-1} \phi^2 dx+\int_{\R^n}\nabla_{\lambda} u \nabla_{\lambda} (\phi^2)|u|^{\gamma-1}u \;dx = \int_{\R^n} |x|^{a}_{\lambda} |u|^{p+\gamma} \phi^2 dx,     $$
therefore
\begin{eqnarray*}
\int_{\R^n} |x|^{a}_{\lambda}  |u|^{p+\gamma} \phi^2 dx&=& \frac{4\gamma}{(\gamma+1)^2} \int_{\R^n}   |\nabla_{\lambda} (|u|^{\frac{\gamma-1}{2}}u)|^2 \phi^2 dx+ \frac{1}{\gamma +1} \int_{\R^n} \nabla_{\lambda} (|u|^{\gamma+1}) \nabla_{\lambda} (\phi^2) dx\nonumber\\&=& \frac{4\gamma}{(\gamma+1)^2} \int_{\R^n}   |\nabla_{\lambda} (|u|^{\frac{\gamma-1}{2}}u)|^2 \phi^2dx-  \frac{1}{\gamma +1} \int_{\R^n} |u|^{\gamma+1} \Delta_{\lambda} (\phi^2) dx .
\end{eqnarray*}
Identity \eqref{e1} then follows by multiplying the latter identity by the factor $\frac{(\gamma +1)^2}{4 \gamma}$.\\\\
\textit{\textbf{Step 2.}} For any $\phi \in C^2_c(\R^n)$ we have
\begin{eqnarray} \label{e2}
\left( p-\frac{(\gamma +1)^2}{4 \gamma}\right)\int_{\R^n} |x|^{a}_{\lambda} |u|^{p+\gamma} \phi^2dx & \leq & \int_{\R^n} |u|^{\gamma+1}\left[ |\nabla_{\lambda} \phi|^2 +\left(\frac{\gamma +1}{4 \gamma}- \frac{1}{2}\right) \Delta_{\lambda} (\phi^2)   \right]dx .
\end{eqnarray}
The function $|u|^{\frac{\gamma-1}{2}}u \phi $ belongs to $C^1_c(\R^n)$, and thus it can be used as a test function in the quadratic form $Q_u$. Hence, the stability assumption on $u$ gives
\begin{eqnarray}\label{e3}
p\int_{\R^n} |x|^{a}_{\lambda} |u|^{p+\gamma} \phi^2 dx \leq \int_{\R^n} |\nabla_{\lambda} (|u|^{\frac{\gamma-1}{2}}u \phi)|^2 dx.
\end{eqnarray}
A direct calculation shows that the right hand side of \eqref{e3} equals to
\begin{multline}\label{e4}
\int_{\R^n} \left[|u|^{\gamma+1} |\nabla_{\lambda} \phi|^2 +  |\nabla_{\lambda}(|u|^{\frac{\gamma-1}{2}}u)|^2 \phi^2 + \frac{1}{2} \nabla_{\lambda} \phi^2 \nabla_{\lambda} (|u|^{\gamma+1}) \right] dx \\=  \int_{\R^n} |u|^{\gamma+1} \left[ |\nabla_{\lambda} \phi|^2 - \frac{1}{2} \Delta_{\lambda} (\phi^2 ) \right] dx + \int_{\R^n} |\nabla_{\lambda}( |u|^{\frac{\gamma-1}{2}}u)|^2 \phi^2 dx.
\end{multline}
From \eqref{e3} and \eqref{e4}, we obtain that
\begin{eqnarray}\label{e5}
p\int_{\R^n} |x|^{a}_{\lambda} |u|^{p+\gamma} \phi^2 dx  &\leq& \int_{\R^n} |u|^{\gamma+1}\left[ |\nabla_{\lambda} \phi|^2 - \frac{1}{2} \Delta_{\lambda} (\phi^2)  \right]dx + \int_{\R^n}  |\nabla_{\lambda} (|u|^{\frac{\gamma-1}{2}}u)|^2 \phi^2 dx.
\end{eqnarray}
Putting this back into \eqref{e1} gives
\begin{eqnarray*}\left( p-\frac{(\gamma +1)^2}{4 \gamma}\right)\int_{\R^n} |x|^{a}_{\lambda} |u|^{p+\gamma} \phi^2dx & \leq & \int_{\R^n} |u|^{\gamma+1}\left[ |\nabla_{\lambda} \phi|^2 +\left(\frac{\gamma +1}{4 \gamma}- \frac{1}{2}\right) \Delta_{\lambda} (\phi^2)  \right]dx .
\end{eqnarray*}
\textit{\textbf{Step 3.}} For any $\gamma \in \left[1,\,\Gamma_M(p)\right)$ and any integer $m\geq \max\{\frac{p+\gamma}{p-1},\; 2\}$ there exists a constant $C(p, m, \gamma)>0$ depending only on $p$, $m$ and $\gamma$
\begin{eqnarray} \label{e6}
\int_{\R^n} |x|^{a}_{\lambda} |u|^{p+\gamma} \psi_R^{2m}dx & \leq & C(p, m, \gamma) \int_{\R^n} |x|^{\frac{-(\gamma+1)a}{p-1}}_{\lambda} \left( |\nabla_{\lambda} \psi_R|^2 + |\psi_R| |\Delta_{\lambda} \psi_R|  \right)^{\frac{p+\gamma}{p-1}}dx ,
\end{eqnarray}
\begin{eqnarray} \label{e07}
\int_{\R^n}  |\nabla_{\lambda} (|u|^{\frac{\gamma-1}{2}}u)|^2 \psi_R^{2m} dx& \leq & C(p, m, \gamma) \int_{\R^n} |x|^{\frac{-(\gamma+1)a}{p-1}}_{\lambda} \left( |\nabla_{\lambda} \psi_R|^2 + |\psi_R| |\Delta_{\lambda} \psi_R|  \right)^{\frac{p+\gamma}{p-1}}dx , \end{eqnarray}
where $\psi_R = \prod_{i=1}^k \psi_{i,R}$. Moreover, the constant $C(p,m,\gamma )$ can be explicitly computed.

From \eqref{e2}, we obtain that
\begin{eqnarray}\label{e7}
\alpha \int_{\R^n} |x|^{a}_{\lambda} |u|^{p+\gamma} \phi^2dx & \leq & \int_{\R^n} |u|^{\gamma+1}|\nabla_{\lambda} \phi|^2 + \beta \int_{\R^n} |u|^{\gamma+1} \Delta_{\lambda} \phi  dx .
\end{eqnarray}
where we have set $\alpha = \left( p-\frac{(\gamma +1)^2}{4 \gamma}\right)$ and $\beta = \frac{1-\gamma}{4 \gamma}$. Notice that $\alpha >0$ and $\beta<0$, since $p>1$ and $\gamma \in \left[1,\,\Gamma_M(p)\right)$.

Now, we set $\phi = \psi_R^m$. The function $\phi$ belongs to $C^2_c(\R^n)$, since $m\geq 2$ and $m$ is an integer, hence it can be used in \eqref{e7}. A direct computation gives
\begin{eqnarray}\label{e8}
\alpha \int_{\R^n} |x|^{a}_{\lambda} |u|^{p+\gamma} \psi_R^{2m} dx & \leq & \int_{\R^n} |u|^{\gamma+1} \psi_R^{2m-2}\left(m^2 |\nabla_{\lambda} \psi_R|^2 + \beta m (m-1) |\nabla_{\lambda} \psi_R|^2+ \beta m \psi_R\Delta_{\lambda} \psi_R  \right)dx, \nonumber\\
\end{eqnarray}
hence
\begin{eqnarray}\label{e9}
\int_{\R^n} |x|^{a}_{\lambda} |u|^{p+\gamma} \psi_R^{2m} dx & \leq & C_1 \int_{\R^n} |u|^{\gamma+1} \psi_R^{2m-2}\left( |\nabla_{\lambda} \psi_R|^2 + |\psi_R ||\Delta_{\lambda} \psi_R|  \right)dx ,
\end{eqnarray}
with $C_1= \frac{m^2+ \beta m (m-1)}{\alpha} > -\frac{\beta m}{\alpha} \geq 0$.

An application of Young’s inequality yields
\begin{eqnarray}\label{e10}
\int_{\R^n} |x|^{a}_{\lambda} |u|^{p+\gamma} \psi_R^{2m} &\leq & C_1 \int_{\R^n} |u|^{\gamma+1} \psi_R^{2m-2}\left( |\nabla_{\lambda} \psi_R|^2 + |\psi_R ||\Delta_{\lambda} \psi_R|  \right)dx \nonumber\\ &=& C_1 \int_{\R^n} |x|^{\frac{(\gamma+1)a}{p+\gamma}}_{\lambda}|u|^{\gamma+1} \psi_R^{2m-2}|x|^{\frac{-(\gamma+1)a}{p+\gamma}}_{\lambda}\left( |\nabla_{\lambda} \psi_R|^2 + |\psi_R ||\Delta_{\lambda} \psi_R|  \right)dx \nonumber\\ & \leq & \frac{\gamma+1}{p+\gamma} \int_{\R^n} |x|^{a}_{\lambda}|u|^{p+\gamma} \psi_R^{(2m-2)\frac{p+\gamma}{\gamma +1}}  + \frac{(p-1) \;C_1}{ p+\gamma} \int_{\R^n}  |x|^{\frac{-(\gamma+1)a}{p-1}}_{\lambda}\left( |\nabla_{\lambda} \psi_R|^2 + |\psi_R ||\Delta_{\lambda} \psi_R|  \right)^{\frac{p+\gamma}{ p-1}}.\nonumber\\
\end{eqnarray}
At this point we notice that $m\geq \max\{\frac{p+\gamma}{p-1},\; 2\}$ implies $(2m-2)\frac{p+\gamma}{p-1}\geq 2m$ and thus
$\psi_R^{(2m-2)\frac{p+\gamma}{\gamma +1}} \leq \psi_R^{2m}$ in $\R^n$, since $0 \leq \psi_R  \leq 1 $ everywhere in $\R^n$.

Therefore, we obtain
\begin{eqnarray*}
\int_{\R^n} |x|^{a}_{\lambda} |u|^{p+\gamma} \psi_R^{2m}  & \leq & \frac{\gamma+1}{p+\gamma} \int_{\R^n} |x|^{a}_{\lambda} |u|^{p+\gamma} \psi_R^{ 2m}  + \frac{(p-1) \;C_1}{ p+\gamma} \int_{\R^n}  |x|^{\frac{-(\gamma+1)a}{p-1}}_{\lambda} \left( |\nabla_{\lambda} \psi_R|^2 + |\psi_R ||\Delta_{\lambda} \psi_R|  \right)^{\frac{p+\gamma}{ p-1}}.\nonumber\\
\end{eqnarray*}
The latter immediately implies
\begin{eqnarray} \label{e12}
\int_{\R^n} |x|^{a}_{\lambda} |u|^{p+\gamma} \psi_R^{2m} dx & \leq & C_1 \int_{\R^n} |x|^{\frac{-(\gamma+1)a}{p-1}}_{\lambda} \left( |\nabla_{\lambda} \psi_R|^2 + |\psi_R ||\Delta_{\lambda} \psi_R|  \right)^{\frac{p+\gamma}{ p-1}}dx,
\end{eqnarray}
which proves inequality \eqref{e6} with $C(p,m,\gamma )= C_1$.

To prove \eqref{e07}, we combine \eqref{e1} and \eqref{e2}. This leads to
\begin{eqnarray*}
\int_{\R^n} |\nabla_{\lambda} (|u|^{\frac{\gamma-1}{2}}u)|^2 \phi^2 dx \leq A \int_{\R^n} |u|^{\gamma +1} |\nabla_{\lambda} \phi|^2 dx + B  \int_{\R^n} |u|^{\gamma +1} \phi \D_{\lambda} \phi dx,
\end{eqnarray*}
where $A= \frac{(\gamma +1)^2}{4 \gamma \alpha }+ \frac{(\gamma +1)}{2 \gamma } >0$ and $B=\frac{\beta (\gamma +1)^2}{4 \gamma \alpha }+ \frac{(\gamma +1)}{2 \gamma } \in \R$.

Now, we insert the test function $\phi = \psi_R^m$ in the latter inequality to find,
\begin{eqnarray*}
\int_{\R^n} |\nabla_{\lambda} (|u|^{\frac{\gamma-1}{2}}u)|^2 \psi_R^{2m} dx \leq  \int_{\R^n} |u|^{\gamma +1} \psi_R^{2m-2} \left(A m^2 |\nabla_{\lambda} \psi_R|^2 + B m (m-1)|\nabla_{\lambda} \psi_R|^2+ B m \psi_R \D_{\lambda} \psi_R \right) dx,
\end{eqnarray*}
and hence
\begin{eqnarray} \label{e11}
\int_{\R^n} |\nabla_{\lambda} (|u|^{\frac{\gamma-1}{2}}u)|^2 \psi_R^{2m}  dx \leq C_2 \int_{\R^n} |u|^{\gamma +1} \psi_R^{2m-2} \left( |\nabla_{\lambda} \psi_R|^2 + |\psi_R| |\D_{\lambda} \psi_R| \right) dx,
\end{eqnarray}
with $C_2= \max\{A m^2 + B m (m-1), \; |B|m \}>0$. \\Using H\"{o}lder’s inequality in \eqref{e11} yields
\begin{eqnarray*}
\int_{\R^n} |\nabla_{\lambda} (|u|^{\frac{\gamma-1}{2}}u)|^2 \psi_R^{2m} &\leq& C_2 \left(\int_{\R^n} |x|^{a}_{\lambda} |u|^{p+\gamma} \psi_R^{(2m-2)\frac{p+\gamma}{\gamma +1}}\right)^{\frac{\gamma+1}{p+\gamma}} \left(\int_{\R^n} |x|^{\frac{-(\gamma+1)a}{p-1}}_{\lambda}\left( |\nabla_{\lambda} \psi_R|^2 + |\psi_R ||\Delta_{\lambda} \psi_R|  \right)^{\frac{p+\gamma}{ p-1}}\right)^{\frac{p-1}{ p+\gamma}}\nonumber\\ &\leq & C_2 \left(\int_{\R^n} |x|^{a}_{\lambda} |u|^{p+\gamma} \psi_R^{2m} \right)^{\frac{\gamma+1}{p+\gamma}} \left(\int_{\R^n} |x|^{\frac{-(\gamma+1)a}{p-1}}_{\lambda} \left( |\nabla_{\lambda} \psi_R|^2 + |\psi_R ||\Delta_{\lambda} \psi_R|   \right)^{\frac{p+\gamma}{ p-1}}\right)^{\frac{p-1}{ p+\gamma}}.
\end{eqnarray*}
Finally, inserting \eqref{e12} into the latter we obtain
\begin{eqnarray*}
\int_{\R^n} |\nabla_{\lambda} (|u|^{\frac{\gamma-1}{2}}u)|^2 \psi_R^{2m} dx &\leq& C_2 C^{\frac{1+\gamma}{ p-1}}_1 \int_{\R^n} |x|^{\frac{-(\gamma+1)a}{p-1}}_{\lambda} \left( |\nabla_{\lambda} \psi_R|^2 + |\psi_R ||\Delta_{\lambda} \psi_R|  \right)^{\frac{p+\gamma}{ p-1}}dx,
\end{eqnarray*}
which gives the desired inequality \eqref{e07}.\\\\
\textit{\textbf{Step 4.}} For any $\gamma \in \left[1,\,\Gamma_M(p)\right)$, there exists a constant $C>0$ independent of $R$ such that
\begin{eqnarray} \label{expr6}
\int_{\Omega_R} \left(|x|^{a}_{\lambda} |u|^{p+\gamma} + |\nabla_{\lambda} (|u|^{\frac{\gamma-1}{2}}u)|^2 \right) dx & \leq & C R^{Q-  \frac{2(p+\gamma)+ (\gamma+1)a }{p-1}},\;\; \forall \; R>0.
\end{eqnarray}

The proof of \eqref{expr6} follows immediately by adding inequality \eqref{e6} to inequality
\eqref{e07} and using Lemma \ref{lem1}.  \qed \\\\
\textbf{Proof of Theorem \ref{th1}.}  By Proposition \ref{proposition1}, there exists a positive constant $C$ independent of $R$ such that
\begin{eqnarray} \label{e13}
\int_{\Omega_R} |x|^{a}_{\lambda} |u|^{p+\gamma} \leq C R^{Q-  \frac{2(p+\gamma)+ a (\gamma+1) }{p-1}}.
\end{eqnarray}
Then it suffices to show that we can always choose a $\gamma \in \left[1,\,\Gamma_M(p)\right)$, such that $Q- \frac{2(p+\gamma)+ a(\gamma+1) }{p-1}<0$. Therefore, by letting $R \rightarrow + \infty$ in \eqref{e13}, we deduce that $$\int_{\R^n} |x|^{a}_{\lambda} |u|^{p+\gamma} =0, $$ which implies that $u \equiv 0$ in $\R^n$.\\

Next, we claim that, under the assumptions on the exponent $p$ assumed in Theorem \ref{th1}, we can always choose $\gamma \in [1,\,\Gamma_M(p))$ such that
 \begin{eqnarray} \label{e09}
 Q-\frac{2(p+\gamma) + a( \gamma +1)}{p-1}<0.
 \end{eqnarray}
As in \cite{Farina}, we consider separately the case $Q\leq 10 + 4a $ and the case $Q> 10 + 4 a$.\\
\textbf{Case 1.} $Q\leq 10 +4 a$ and $p>1$. In this case we have $$2(p+\Gamma_M(p)) + a( \Gamma_M(p) +1)> 2(3p-1+2(p-1)) + a(2p + 2(p-1) >(10+ 4 a)(p-1)$$
and therefore
\begin{eqnarray} \label{e010}
Q-\frac{2(p+\Gamma_M(p)) + a( \Gamma_M(p) +1)}{p-1}<Q - (10 + 4 a)\leq0.
\end{eqnarray}
The latter inequality and the continuity of the function $x\mapsto Q-\frac{2(p+x)+a(x+1)}{p-1}$ immediately imply the existence of $\gamma \in [1,\,\Gamma_M(p))$ satisfying \eqref{e09}.\\\\
\textbf{Case 2.} $Q> 10 + 4 a$ and $1<p<p_c(Q,a)$. In this case we consider the real-valued function $x\mapsto g(x):= \frac{2(x+\Gamma_M(x)) + a(\Gamma(x)+1)}{x-1}$ on $(1,\,+\infty)$. Since $g$ is strictly decreasing function satisfying $\lim_{x\rightarrow 1^+}g(x)=+\infty$ and $\lim_{x\rightarrow +\infty}g(x)=10+4a$, there exists a unique $p_0>1$ such that $Q=g(p_0) $. We claim that $p_0=p_c(Q,a)$. Indeed,
$$Q=g(p)\;\Leftrightarrow\; (Q-2)(p-1)- (4+2a)p=(4+2a)\sqrt{p(p-1)}$$\;$$\Leftrightarrow\; (Q-10-4a)(Q-2)p^2+(-2(Q-2)^2+ 4(a+2)(Q+a))p+(Q-2)^2= 0, $$
which implies that
\begin{eqnarray} \label{e011}
(Q-10-4a)(Q-2)p_0^2+(-2(Q-2)^2+ 4(a+2)(Q+a))p_0+(Q-2)^2= 0,
\end{eqnarray}
 and
\begin{eqnarray} \label{e012}
(Q-2)(p_0-1)- (4+2a)p_0 >(4+2a)(p_0-1).
\end{eqnarray}
The roots of \eqref{e011}
\begin{eqnarray} \label{e013}
p_1=\frac{(Q-2)^2-2(a+2)(a+Q)+2 \sqrt{(a+2)^3(a+2Q-2)}}{(Q-2)(Q-4a-10)}= p_c(Q,a),
\end{eqnarray}
\begin{eqnarray} \label{e014}
p_2=\frac{(Q-2)^2-2(a+2)(a+Q)-2 \sqrt{(a+2)^3(a+2Q-2)}}{(Q-2)(Q-4a-10)}< p_0,
\end{eqnarray}
while \eqref{e012} easily implies $p_0>\frac{Q-6-2a}{Q-4a-10}>p_2 $. This proves that $p_0= p_1$. Hence $$p_c(Q,a)=\frac{(Q-2)^2-2(a+2)(a+Q)+2 \sqrt{(a+2)^3(a+2Q-2)}}{(Q-2)(Q-4a-10)}$$ as claimed. Since we have just proven that $g(p_c(Q,a))=Q$ and $g$ is a strictly decreasing function, it follows that
 \begin{eqnarray} \label{e015}
 \forall \; 1<p<p_c(Q,a), \; Q<g(p).
 \end{eqnarray}
 Now we can conclude as in the first case, i.e, the continuity of $x\mapsto Q-\frac{2(p+x) + a(x+1)}{p-1}$ immediately implies the existence of $\gamma \in [1,\,\Gamma_M(p))$ satisfying \eqref{e09}. \qed
\section{\large{The Liouville theorem for solutions which are stable outside a compact set of $\R^n$: proof of Theorem \ref{th2}}}

  In this section, we prove Proposition \ref{proposition2} and Theorem \ref{th2}. \\\\
\textbf{Proof of Proposition \ref{proposition2}.} Let $u \in C^2(\R^n)$ be a solution of \eqref{equation} and $\phi \in C^1_c(\Omega_R)$. Multiplying equation \eqref{equation} by $T(u) \phi $ and integrating by parts in $\Omega_R$, we obtain
\begin{eqnarray} \label{p3}
-\int_{\Omega_R} \D_{\lambda} u T(u) \phi dx &=&-\int_{\Omega_R} \D_{\lambda} u \;\epsilon_j\; x^{(j)}\; \nabla_{x^{(j)}} u \;\phi dx \nonumber\\&=& \int_{\Omega_R} \lambda^2_{i} \;\nabla_{x^{(i)}} u\; \nabla_{x^{(i)}} \left( \epsilon_j x^{(j)} \nabla_{x^{(j)}} u \phi \right) dx \nonumber\\&=& \int_{\Omega_R} \lambda^2_{i} \;\nabla_{x^{(i)}} u \;\epsilon_j \delta_{ij} \;\nabla_{x^{(j)}} u \; \phi dx  + \int_{\Omega_R} \lambda^2_{i} \; \nabla_{x^{(i)}} u \; \epsilon_j\;  x^{(j)} \; \nabla_{x^{(i)}}\left(\nabla_{x^{(j)}} u \right) \; \phi dx \nonumber\\&&+\;\; \int_{\Omega_R} \lambda^2_{i} \; \nabla_{x^{(i)}} u \; \epsilon_j \; x^{(j)} \; \nabla_{x^{(j)}} u \; \nabla_{x^{(i)}}\phi dx\nonumber\\&:=& I_1 + I_2 + I_3,
\end{eqnarray}
Here and in the sequel, we use the Einstein summation convention: an index occurring twice in a product is to be summed from $1$ up to the space dimension.

Obviously
\begin{eqnarray} \label{p4}
I_1 :&=& \int_{\Omega_R} \lambda^2_{i} \nabla_{x^{(i)}} u \; \epsilon_j  \delta_{ij}\;  \nabla_{x^{(j)}} u \; \phi dx \nonumber\\&=& \int_{\Omega_R} \lambda^2_{i} \left|\nabla_{x^{(i)}} u \right|^2 \epsilon_i  \phi dx.
\end{eqnarray}
Moreover, an integration by parts in $I_2$ gives
\begin{eqnarray*}
I_2 :&=& \int_{\Omega_R} \lambda^2_{i} \nabla_{x^{(i)}} u \epsilon_j x^{(j)} \nabla_{x^{(i)}}\left(\nabla_{x^{(j)}} u \right) \phi dx \nonumber\\&=& - \int_{\Omega_R} \nabla_{x^{(j)}} (\lambda^2_{i}) \left|\nabla_{x^{(i)}} u \right|^2 \epsilon_j x^{(j)} \phi dx -I_2- \int_{\Omega_R} \lambda^2_{i} \left|\nabla_{x^{(i)}} u \right|^2 \epsilon_j n_j \phi dx -\int_{\Omega_R} \lambda^2_{i} \left|\nabla_{x^{(i)}} u \right|^2 \epsilon_j x^{(j)} \nabla_{x^{(j)}}\phi dx \nonumber\\&=& - 2\int_{\Omega_R} \lambda_{i} \left|\nabla_{x^{(i)}} u \right|^2 T(\lambda_{i}) \phi dx -I_2- Q \int_{\Omega_R} |\nabla_{\lambda}u |^2 \phi dx  -\int_{\Omega_R} |\nabla_{\lambda}u |^2 T(\phi)dx.
\end{eqnarray*}
Since $\lambda_i$ is $\delta_t$-homogeneous of degree $\epsilon_i-1$, then $T(\lambda_{i}) = (\epsilon_i-1) \lambda_{i} $. Hence
\begin{eqnarray*}
I_2&=& - 2\int_{\Omega_R} (\epsilon_i-1)\lambda^2_{i} \left|\nabla_{x^{(i)}} u \right|^2 \phi dx -I_2- Q \int_{\Omega_R} |\nabla_{\lambda}u |^2 \phi dx-\int_{\Omega_R} |\nabla_{\lambda}u |^2 T(\phi)dx \nonumber\\&=&
(2-Q) \int_{\Omega_R} |\nabla_{\lambda}u |^2 \phi dx -2I_1-I_2- \int_{\Omega_R} |\nabla_{\lambda}u |^2 T(\phi)dx.
\end{eqnarray*}
Then
\begin{eqnarray} \label{p5}
I_2= \frac {2-Q}{2} \int_{\Omega_R} |\nabla_{\lambda}u |^2 \phi dx -I_1 - \frac 12 \int_{\Omega_R} |\nabla_{\lambda}u |^2 T(\phi)dx.
\end{eqnarray}
It is easily seen that
\begin{eqnarray} \label{p6}
I_3:&=&\int_{\Omega_R} \lambda^2_{i} \nabla_{x^{(i)}} u \epsilon_j x^{(j)} \nabla_{x^{(j)}} u \nabla_{x^{(i)}} \phi dx\nonumber\\&=& \int_{\Omega_R} \nabla_{\lambda}u \nabla_{\lambda} \phi  T(u) dx.
\end{eqnarray}
Hence, by \eqref{p3},
\begin{eqnarray} \label{point1}
-\int_{\Omega_R} \D_{\lambda} u T(u) \phi dx = \frac{2- Q}{2} \int_{\Omega_R} |\nabla_{\lambda}u |^2 \phi dx - \frac 12 \int_{\Omega_R} |\nabla_{\lambda}u |^2 T(\phi) dx+ \int_{\Omega_R} \nabla_{\lambda}u \nabla_{\lambda} \phi T(u)dx.\nonumber\\
\end{eqnarray}
On the other hand, an integration by parts gives
\begin{multline*}
\int_{\Omega_R} |x|^{a}_{\lambda} |u|^{p-1} u T(u) \phi dx = \frac{1}{p+1} \int_{\Omega_{R}} |x|^{a}_{\lambda} \nabla_{x^{(j)}}(|u|^{p+1}) \epsilon_j x^{(j)} \phi dx \\= -\frac{Q}{p+1} \int_{\Omega_{R}} |x|^{a}_{\lambda} |u|^{p+1} \phi -\frac{a}{p+1} \int_{\Omega_R} |x|^{a-1}_{\lambda} |u|^{p+1}T(|x|_{\lambda})\phi-\frac{1}{p+1} \int_{\Omega_R} |x|^{a}_{\lambda} |u|^{p+1} T(\phi) dx.
\end{multline*}
If $T(|x|_{\lambda})= |x|_{\lambda}$, then
\begin{eqnarray} \label{point2}
\int_{\Omega_R} |x|^{a}_{\lambda} |u|^{p-1} u T(u) \phi dx &=& \frac{1}{p+1} \int_{\Omega_R} |x|^{a}_{\lambda} \nabla_{x^{(j)}}(|u|^{p+1}) \epsilon_j x^{(j)} \phi dx \nonumber\\&=& -\frac{Q+a}{p+1} \int_{\Omega_R} |x|^{a}_{\lambda} |u|^{p+1} \phi -\frac{1}{p+1} \int_{\Omega_R} |x|^{a}_{\lambda} |u|^{p+1} T(\phi) dx.
\end{eqnarray}
Clearly \eqref{pohozaev} follows directly from \eqref{point1} and \eqref{point2}. \qed \\\\
\textbf{Proof of Theorem \ref{th2}.}  Let $u$ be a solution of \eqref{equation} which is stable outside a compact set. We begin defining some smooth compactly supported functions which will be used several times in the sequel. More precisely, for $R_*>0$, we choose a function $\zeta_{i,R}\in C^2_c(\R^{n_i})$, $i=1,...,k$, $0\leq \zeta_{i,R} \leq 1$, everywhere on $\R^{n_i}$ and
$$\begin{cases}
\zeta_{i,R}(x^{(i)})=0 \;\;&\text{if \;\; $|x^{(i)}|<  R_*+1$ \; or \;$|x^{(i)}|> 2 R^{\epsilon_i} $}, \\ \zeta_{i,R}(x^{(i)})=1 \;\;\; &\text{if\;\;$R_* +2<|x^{(i)}|< R^{\epsilon_i} $}, \\
|\nabla_{x^{(i)}} \zeta_{i,R}|^2 + |\D_{x^{(i)}} \zeta_{i,R}|\leq CR^{-2 \epsilon_i}\;\; &\text{for \;\; $\{R^{\epsilon_i} < |x^{(i)}| < 2 R^{\epsilon_i}\}$}.
\end{cases}$$

The rest of the proof splits into several steps.\\
\textbf{\textit{Step 1.}} Let $p>1$. There exists $R_*>0$ such that for every $\gamma \in \left[1,\,\Gamma_M(p)\right)$ and every $R^{\epsilon_i}>R_* +2$, we have \begin{eqnarray} \label{e14}
\int_{\Sigma_0(R)} \left(|x|^{a}_{\lambda} |u|^{p+\gamma}+ |\nabla_{\lambda} ( |u|^{\frac{\gamma-1}{2}}u)|^2 \right)  dx \leq C_{R_*}+C R^{Q- \frac{2(p+\gamma) + (\gamma +1) a}{p-1}},
\end{eqnarray}
where $\Sigma_0(R)= \Omega_R \backslash B_1(0, R_*+2)\times ...\times B_k(0, R_*+2)$, $C_{R_*}$ and $C$ are positive constants depending on $p$, $\gamma $, $R_*$ but not on $R$.


Since $u$ is stable outside a compact set of $\R^n$, there exists $R_*>0$ such that, similar to that of Proposition \ref{proposition1} we derive
\begin{eqnarray*}
\int_{\Sigma_0(R)} \left(|x|^{a}_{\lambda} |u|^{p+\gamma}+ |\nabla_{\lambda} ( |u|^{\frac{\gamma-1}{2}}u)|^2 \right) dx &\leq& C(p,m,\gamma)\int_{\R^n} |x|^{\frac{-a (\gamma +1)}{p-1}}_{\lambda}\left(|\nabla_{\lambda} \zeta_{R}|^2 + |\zeta_{R} | |\D_{\lambda} \zeta_{R} |\right)^{\frac{p+\gamma}{p-1}} dx \nonumber\\&\leq& C_{R_*} + C R^{Q-  \frac{2(p+\gamma)+ (\gamma +1)a}{p-1}},
\end{eqnarray*}
where $\zeta_{R}=\prod_{i=1}^n \zeta_{i, R}$. Hence, the desired integral estimate \eqref{e14} follows.\\\\
\textbf{\textit{Step 2.}} If  $Q=2$ and $1<p<+\infty$ or $Q\geq 3$ and $1<p<\frac{Q+2+ 2a }{Q-2}$, then $u\equiv 0$.\\

By choosing $\gamma=1$  and using \textit{Step $1$}, we get $|x|^{\frac{a}{p+1}}_{\lambda} u\in L^{p+1}(\R^n)$ and $|\nabla_{\lambda} u|\in L^{2}(\R^n)$ for $1<p<p_s(Q,a)$.\\

Take $\phi = \psi_{ R}= \prod_{i=1}^k \psi_{i,R}$ in \eqref{pohozaev} where  $\psi_{i,R}$ defined as above. Since $|x|^{\frac{a}{p+1}}_{\lambda} u\in L^{p+1}(\R^n)$ and $|\nabla_{\lambda} u|\in L^{2}(\R^n)$, then
 \begin{eqnarray} \label{homogene} \int_{\Sigma_R}|\nabla_{\lambda} u|^2 dx \rightarrow 0, \; \mbox{as}\; R \rightarrow +\infty \quad \mbox{ and } \quad \int_{\Sigma_R} |x|^{a}_{\lambda} |u|^{p+1}dx \rightarrow 0, \; \mbox{as}\; R \rightarrow +\infty,\end{eqnarray}
where $\Sigma_R= \Omega_{2R} \backslash \Omega_{R}$.\\
Recalling that $\lambda_i$ and $\lambda_i \nabla_{x^{(i)}} u$ are $\delta_t$-homogeneous of degree $\epsilon_i-1$ and one respectively. Then, since $T$ generates $(\delta_t)_{t\geq 0}$, we have
\begin{eqnarray} \label{Euler}
T(\lambda_i)= (\epsilon_i-1) \lambda_i \quad \mbox{and} \quad T(\lambda_i \nabla_{x^{(i)}} u)= \lambda_i \nabla_{x^{(i)}} u.
\end{eqnarray}
 Integrating by parts and using \eqref{Euler}, we derive
\begin{eqnarray} \label{integration}
 &&\int_{\Omega_{2R}} \nabla_{\lambda} u \nabla_{\lambda} \psi_{ R}  T(u) = \int_{\Omega_{2R}}  \lambda_i \nabla_{x^{(i)}} u \lambda_i \nabla_{x^{(i)}} \psi_{ R} \epsilon_j x^{(j)} \nabla_{x^{(j)}} u\nonumber\\ &=&  - \int_{\Omega_{2R}}  T(\lambda_i \nabla_{x^{(i)}} u) \lambda_i \nabla_{x^{(i)}} \psi_{ R}\; u  - \int_{\Omega_{2R}} \lambda_i \nabla_{x^{(i)}} u T( \lambda_i) \nabla_{x^{(i)}} \psi_{ R} \;u  -\int_{\Omega_{2R}}  \lambda_i^2 \nabla_{x^{(i)}} u T(\nabla_{x^{(i)}} \psi_{ R})  \; u\nonumber\\&& -\;\;  Q \int_{\Omega_{2R}}  \nabla_{\lambda} u \nabla_{\lambda}\psi_{ R}  \; u \nonumber\\&=& -(Q+1) \int_{\Omega_{2R}}  \nabla_{\lambda} u \nabla_{\lambda}\psi_{ R}  \; u - \int_{\Omega_{2R}} (\epsilon_i-1)\lambda_i^2 \nabla_{x^{(i)}} u  \nabla_{x^{(i)}} \psi_{ R} \;u -\int_{\Omega_{2R}}  \lambda_i^2 \nabla_{x^{(i)}} u T(\nabla_{x^{(i)}} \psi_{ R})  \; u \nonumber\\&=&  \frac{Q+1}{2} \int_{\Omega_{2R}}  u^2 \D_{\lambda}\psi_{ R} + \int_{\Omega_{2R}} \frac{\epsilon_i-1}{2}  u^2  \lambda_i^2 \D_{x^{(i)}} \psi_{ R} + \frac 12 \int_{\Omega_{2R}}   u^2 \lambda_i^2 \nabla_{x^{(i)}}[T(\nabla_{x^{(i)}} \psi_{ R})]
\end{eqnarray}
By Lemma \ref{lem1}, \eqref{integration} and using H\"{o}lder's inequality, we obtain
\begin{eqnarray}\label{label1}
\left| \int_{\Omega_{2R}} \nabla_{\lambda} u\nabla_{\lambda} \psi_{ R}  T(u)\right| &\leq& \frac{C}{R^{-2}} \int_{\Sigma_R} u^2 \nonumber\\&=& \frac{C}{R^{-2}} \int_{\Sigma_R} |x|^{\frac{-2a}{p+1}}_{\lambda} |x|^{\frac{2a}{p+1}}_{\lambda}  u^2 \nonumber\\&\leq& C R^{(Q-\frac{2a}{p-1}) \frac{p-1}{p+1}-2} \left(\int_{\Sigma_R}|x|^{a}_{\lambda} |u|^{p+1}\right)^{\frac{2}{p+1}}.
\end{eqnarray}
Similarly, we get
\begin{eqnarray} \label{label2}
\left| \int_{\Omega_{2R}}
\left[- \frac 12 |\nabla_{\lambda} u|^2 + \frac{|x|^{a}_{\lambda}}{p+1} |u|^{p+1}\right] T(\psi_{ R})  \right| \leq C \int_{\Sigma_R}\left(|\nabla_{\lambda} u|^2 + |x|^{a}_{\lambda} |u|^{p+1}\right).
\end{eqnarray}
From \eqref{homogene}, \eqref{label1} and \eqref{label2}, we obtain
\begin{eqnarray*} \label{e016}
\lim_{R \rightarrow +\infty}\left| \int_{\Omega_{2R}} \left( \nabla_{\lambda} u \nabla_{\lambda} \psi_{ R}  T(u) +\left[- \frac 12 |\nabla_{\lambda} u|^2 + \frac{|x|^{a}_{\lambda}}{p+1} |u|^{p+1}\right] T(\psi_{ R}) \right) \right| = 0.
\end{eqnarray*}
As a consequence, \eqref{pohozaev} becomes
\begin{eqnarray} \label{e17}
 \frac{Q-2}{2} \int_{\R^n} |\nabla_{\lambda} u|^2 dx -\frac{Q+a}{p+1} \int_{\R^n} |x|^{a}_{\lambda} |u|^{p+1}  dx = 0.
\end{eqnarray}
On the other hand, multiplying equation \eqref{equation} by $u \psi_{ R} $ and integrating by parts yields
\begin{eqnarray*}
\int_{\R^n} |\nabla_{\lambda} u|^2 \psi_{ R} dx - \int_{\R^n} |x|^{a}_{\lambda} |u|^{p+1}  \psi_{ R} dx = \frac{1}{2} \int_{\R^n} u^2 \D_{\lambda} \psi_{ R}dx.
\end{eqnarray*}
Since $1<p<p_s(Q,a)$, we get
\begin{eqnarray*}
\left|\int_{\R^n} u^2 \D_{\lambda} \psi_{ R}dx\right| &\leq& \left(\int_{\R^n} |x|^{a}_{\lambda} |u|^{p+1}dx \right)^{\frac{2}{p+1}} \left(\int_{\Sigma_R} |x|^{\frac{-2a}{p-1}}_{\lambda}|\D_{\lambda} \psi_{ R}|^{\frac{p+1}{p-1}}dx \right)^{\frac{p-1}{p+1}} \nonumber\\ &\leq& C R^{Q\frac{p-1}{p+1}- 2 -\frac{2 a}{p+1}} \rightarrow 0 \; \mbox{as}\; R \rightarrow +\infty.
\end{eqnarray*}
Then

\begin{eqnarray} \label{e18}
\int_{\R^n} |\nabla_{\lambda} u|^2 dx= \int_{\R^n} |x|^{a}_{\lambda} |u|^{p+1}dx .
\end{eqnarray}
To complete the proof we combine \eqref{e17} and \eqref{e18} to get $$\left(\frac{Q-2}{2}- \frac{Q+a}{p+1}\right)\int_{\R^n} |x|^{a}_{\lambda} |u|^{p+1}dx=0, $$
but $\frac{Q-2}{2} - \frac{Q+a}{p+1} \neq0$, since $p$ is subcritical, hence $u$ must be identically zero, as claimed. \qed

\end{large}

\section*{\large{References}}

 \begin{thebibliography}{777}
 \bibliographystyle{alpha}
\begin{large}
\bibitem{AnhMy} C. T. Anh and B. K. My,\emph{Liouville-type theorems for elliptic inequalities involving the
$\D_\lambda$-Laplace operator,} Complex Variables and Elliptic Equations, (2016) http://dx.doi.org/10.1080/17476933.2015.1131685.

\bibitem{Bozhkov} Y.Bozhkov and E.Mitidieri,\emph{Conformal Killing vector fields and Rellich type identies on Riemannian manifollds, } I, in: Geometric Methods in PDE’s, 65–80, in: Lect. Notes Semin. Interdiscip. Mat., Potenza, \textbf{7} (2008) 65–80.

\bibitem{Cauchy} A.Cauchy, \emph{M\'{e}moires sur les fonctions compl\'{e}mentaires [Memoirs on complementary functions], }
C. R. Acad. Sci. Paris, \textbf{19} (1844) 1377–1384.

\bibitem{Chen} WX.Chen and C.Li, \emph{Classification of solutions of some nonlinear elliptic equations}. Duke Math. J.
\textbf{63} (1991) 615–622.

\bibitem{Cuccu} F.Cuccu, A.Mohammed and G.Porru, \emph{ Extensions of a theorem of Cauchy–Liouville}. J. Math. Anal. Appl. \textbf{369} (2010) 222–231.

\bibitem{Dolcetta} IC.Dolcetta and A.Cutrì, \emph{On the Liouville property for the sublaplacians}. Ann. Scuola Norm. Sup.
Pisa Cl. Sci.  \textbf{25} (1997) 239–256.

\bibitem{Farina}
A.Farina, \emph{On the classification of solutions of the Lane–Emden equation on unbounded domains of $\R^n$,} J. Math. Pures Appl. \textbf{87} (2007) 537–561.

\bibitem{Franchi} B.Franchi and E.Lanconelli, \emph{Une m\'{e}trique associ\'{e}e à une classe d'op\'{e}rateurs elliptiques d\'{e}g\'{e}n\'{e}r\'{e}s,
(French) [A metric associated with a class of degenerate elliptic operators]}. Conference on linear partial and pseudodifferential operators (Torino, 1982). Rend. Sem. Mat. Univ. Politec. Torino
1983, Special Issue. (1984) (1983) 105–114.

\bibitem{Garofalo} N.Garofalo and E.Lanconelli, \emph{Existence and nonexistence results for semilinear differential equations on the Heisenberg group}. Indiana Univ. Math. J. \textbf{41} (1992) 71–98.

\bibitem{Gidas1} B.Gidas and J.Spruck, \emph{Global and local behavior of positive solutions of nonlinear elliptic equations}. Commun. Pure Appl. Math. \textbf{34} (1981) 525–598.

\bibitem{Gidas2} B.Gidas and J.Spruck, \emph{A priori bounds for positive solutions of a nonlinear elliptic equations}.
Commun. Partial Differ. Equ. \textbf{6} (1981) 883–901.

\bibitem{Kogoj3} A.E.Kogoj and E.Lanconelli, \emph{On semilinear $\D_{\lambda}$-Laplace equation}. Nonlinear Anal. \textbf{75} (2012) 4637–4649.

\bibitem{Kogoj2} A.E.Kogoj and E.Lanconelli, \emph{Liouville theorem for $X$-elliptic operators}. Nonlinear Anal.
\textbf{70} (2009) 2974–2985.

\bibitem{Kogoj1} A.E.Kogoj and E.Lanconelli, \emph{ Liouville theorems in half-spaces for parabolic hypoelliptic equations}. Ric. Mat. \textbf{55} (2006) 267–282.


\bibitem{Liouville} J.Liouville, C.R.Acad. Sci. Paris. \textbf{19} (1844) 1262.

\bibitem{Monticelli} DD.Monticelli, \emph{Maximum principles and the method of moving planes for a class of degenerate elliptic linear operators}. J. Eur. Math. Soc. \textbf{12} (2010) 611–654.

\bibitem{Monti} R.Monti and D.Morbidelli, \emph{Kelvin transform forGrushin operators and critical semilinear equations}.
Duke Math. J. \textbf{131} (2006) 167–202.

\bibitem{Pohozaev} S.I.Pohozaev, \emph{Eigenfunctions of the equation $\D u + \lambda f (u) = 0$}. Dokl. Akad. Nauk SSSR \textbf{165} (1) (1965) 33–36.

\bibitem{Polacik} P.Pol\`{a}\v{c}ik, P.Quittner and P.Souplet,   \emph{ Singularity and decay estimates in superlinear problems via Liouville-type theorems, I: elliptic equations and systems}. Duke Math. J. \textbf{139} (2007) 555–579.

\bibitem{Pucci} P.Pucci and J.Serrin, \emph{A general variational identity}. Indiana Univ. Math. J. \textbf{35} (3) (1986) 681–703.

\bibitem{Quittner} P.Quittner and P.Souplet, \emph{ Superlinear parablolic problems: blow-up, global existence and steady
states}. Basel: Birkhäuser Verlag; 2007.

\bibitem{Serrin} J.Serrin and H.Zou, \emph{ Cauchy–Liouville and universal boundedness theorems for quasilinear elliptic
equations and inequalities}. Acta Math. \textbf{189} (2002) 79–142.

\bibitem{WY} C.Wang and D.Ye, \emph{Some Liouville theorems for H\'{e}non type elliptic equations,} Journal of Functional Analysis \textbf{262} (2012) 1705-1727.

\bibitem{Yu} X.Yu, \emph{Liouville type theorem for nonlinear elliptic equation involving Grushin operators}. Commun. Contem. Math. \textbf{17} (2015) 1450050 (12$p$).

\end{large}
\end{thebibliography}
\end{document}